\documentclass[11pt]{article}

\RequirePackage[colorlinks]{hyperref}
\usepackage{graphicx}
\usepackage{epsfig}
\usepackage{amsmath}
\usepackage{amsthm}
\usepackage{amssymb}
\usepackage{verbatim}
\usepackage{url}

\allowdisplaybreaks

\newtheorem{thm}{Theorem}[section]

\newtheorem{definition}{Definition}
\newtheorem{proposition}[thm]{Proposition}
\newtheorem{theorem}{Theorem}[section]
\newtheorem{lemma}[thm]{Lemma}

\usepackage{verbatim}
\title{Optimal cash management using impulse control}
\author{~\\
  {\rm
    \begin{tabular}{c}
      Peter Lakner and Josh Reed
      \\
      Leonard N. Stern School of Business  \\
      New York University  \\
     44 West 4th St. \\
      New York, NY 10012 \\
      ~~\\
    plakner@stern.nyu.edu,~jreed@stern.nyu.edu
    \end{tabular}
  }}

\date{\today}

\begin{document}

\maketitle

\begin{abstract}We consider the impulse control of Levy processes under the infinite horizon, discounted
cost criterion. Our motivating example is
the cash management problem in which a controller is charged a fixed plus proportional cost
for adding to or withdrawing from his/her reserve, plus an opportunity cost for keeping any
cash on hand. Our main result is to
provide a verification theorem for the optimality of control band policies in this scenario. We also
analyze the transient and steady-state behavior of the controlled process under control band policies and explicitly solve for the optimal policy in the case in which the Levy process to be controlled is
the sum of a Brownian motion with drift and a compound Poisson process with exponentially distributed jump sizes.
\end{abstract}

\section{Introduction}

Impulse control problems have a long history related to applications to the cash management problem.  In the present paper, we consider the impulse
control of Levy processes. Our motivating application is the cash
management problem in which there exists a system manager who must control the amount of cash he/she has on hand. We assume that the manager's cash on hand fluctuates due to randomly occurring withdrawals from and deposits to his/her account but
that the manager is charged a fixed plus proportional cost for any specific, intentional adding to or withdrawing from his/her reserves and that
there exists an opportunity cost for keeping too little or too much cash on hand. The manager's objective is to minimize his/her long run
opportunity cost of keeping cash on hand plus any cost incurred from depositing or withdrawing from
the reserve. An alternative motivating application which is also considered in the literature is a manager
who wishes to control his/her inventory level. The manager's inventory level fluctuates randomly and he/she may increase or
decrease his/her inventory level at will by expediting or salvaging parts, paying a fixed plus proportional cost to do so. The
manager's objective is to minimize his/her long run inventory holding costs plus costs of expediting and salvaging.

Our first main result in the paper is to provide a verification theorem for the
optimality of control band policies for the impulse control of Levy processes. Our result is fairly general
and holds for a wide class of opportunity cost functions. We then explicitly calculate the Laplace transform with respect to time and steady-state distribution of any spectrally positive Levy process controlled under a control band policy. In Section \ref{Sec::Example}, we consider the special case of a Levy process which is comprised of the sum of a Brownain motion and a compound Poisson process with exponentially distributed jump sizes.
In this specific case, we show how one may use the results derived in this paper in order to characterize the value function for the associated control problem, and characterize the band levels in the optimal control policy as a solution of a system of equations. Moreover, we also show  how one may determine the steady-state distribution of the controlled process when the underlying Levy process is the sum of a negative drift
and a compound Poisson process with exponentially distributed jump sizes.

The technique of impulse control was originally developed by Bensoussan and Lions \cite{BL73,BL75} and extended
by Richard et. al \cite{CR78,richard1976optimal}. Harrison, Selke and Taylor \cite{HST83}  and Sulem \cite{sulem1986solvable} also consider the impulse
control of Brownian motion and  explicitly calculate the critical parameters determining the
optimal policy. In \cite{feng2010computational}, an iterative computational scheme is to provided
in order to determine the optimal policy for the impulse control of Brownain motion. Recently, Ormeci, Dai and Vande Vate \cite{ODV08} have considered impulse control of Brownian motion under
the average cost criterion and again show that a control band policy is optimal. Cadenillas, Zapatero, and Sarkar \cite{CZS} and Cadenillas, Lakner, and Pinedo \cite{CLP} solved the Brownian case with a mean-reverting drift.
None of the above mentioned works allow there to be jumps in the process to be controlled. In \cite{bensoussan2006optimality}, the optimality of an $(s,S)$ policy is proven for a process which is the
sum of a constant drift, a Brownian motion, and a compound Poisson process.
In related work, Bar-Ilan, Perry and Stadje \cite{BPS04} have also considered
the problem of impulse control of Levy processes for the specific case in which the Levy process is a sum of a Brownian motion and a compound Poisson process. Assuming
that a control band policy is optimal, their main results are to evaluate the cost functionals of the resulting policy through a fundamental identity derived from the martingale originally introduced by Kella and Whitt \cite{KW92}. In \cite{yamazaki2017inventory}, Yamazaki uses a scale function approach to study one-sided impulse control problems for spectrally positive Levy demand processes. In \cite{czarna2020optimality}, an impulse control problem is studied for a refracted Levy process where the ruin time is modeled by a Parisian delay. Finally, \cite{wu2019optimal} studies the impulse control of a geometric Levy process.

\section{The Model}
\label{Sec::Model}

In this Section we provide the specifics of the model described in the Introduction.   All forthcoming processes are  assumed to live on a probability triplet equipped with a filtration ${\cal F}=\{{\cal F}_t,\ 0\le t<\infty\}$. We begin by assuming that $Y_t$ is a Levy process started from $x$
with Levy measure $\nu$ such that
\begin{equation}\int_{\left\{|y|\ge 1\right\}}|y|\nu(dy)<\infty. \label{Main}
\end{equation}
The process $Y_t$ will be used to represent the cash on hand process assuming that the manager exerts no control
by making no deposits to or withdrawal from his/her fund.
Let $J(\omega,dt,dy)=J(dt,dy)$ be the jump measure of $Y$. Then, $Y_t$ has the Ito-Levy  decomposition
$$Y_t=x+\mu t+\sigma w_t +A_t+M_t$$
where $M_t$ is the following martingale
$$M_t=\int_{\left\{|x|<1\right\}}x \left\{J((0,t],dx)-t\nu(dx)\right\}$$
and $A$ is the sum of the \lq\lq large\rq\rq\ jumps
$$A_t=\sum_{0<s\le t}\Delta Y_s1_{\left\{|\Delta Y_s|\ge 1\right\}}.$$
The process $w$ is assumed to be a standard Wiener process and $\mu$ is a constant. We do not make any assumption regarding $\sigma$, it may be zero or non-zero. Also we allow $\nu(\Re)=0$ in which case $Y$ is continuous. The case when both $\sigma=0$ and $\nu(\Re)=0$ is also included, although this case is trivial ($Y$ is deterministic in this case). In general, we will use $P_x$ to denote the probability measure under which $Y_t$ is started  from $x$ and $E_x$ its associated expectation operator.
We note that $M$ is a quadratic pure jump local martingale (\cite{protter}, Chapter II, Section 6) with quadratic variation
$$\left[M\right]_t=\sum_{0<s\le t}\left(\Delta Y_s\right)^21_{\left\{|\Delta Y_s|\le 1\right\}}=\int_{(0,t]\times\Re}y^21_{\{|y|<1\}}J(ds,dy),$$
which has expected value
\begin{equation}E\left[M\right]_t= t\int_\Re y^21_{\{|y|<1\}}\nu(dy)<\infty.\label{new7} \end{equation}
It follows (\cite{protter}, Chapter II, Corollary 3 to Theorem 27) that $M$ is a square-integrable martingale.

We let $$(T,\Xi)=(\tau_1,\tau_2,\dots,\tau_n,\dots,\ \xi_1,\xi_2,\dots,\xi_n,\dots,)$$
denote the impulse control policy used by the manager where
$0\le\tau_1<\tau_2<\tau_3\dots$ are stopping times and $\xi_n\not=0$\ (a.s.) is an ${\cal F}_{\tau_n}$ measurable random variable for each $n\ge 1$. Positive values of $\xi_n$ represent deposits by the manager into his/her fund
and negative values represent withdrawal. In order to simplify future developments we define  $\tau_0=0$ and $\xi_0=0$. In principle we also allow that with positive probability there are only finitely many interventions. However, we introduce the following convention which essentially changes nothing, but simplifies future formulas. Let $\Lambda(\omega)\le\infty$ be the number of interventions. For those $\omega$'s such that $\Lambda(\omega)<\infty$ we define
\begin{equation}\tau_{\Lambda(\omega)+k}=\infty \quad\hbox{and}\quad \xi_{\Lambda(\omega)+k}=0\label{new1}
\end{equation}
for $k\ge 1$.

As described in the Introduction, the controlled cash on hand process $X_t$ follows the dynamics
$$X_t=Y_t+ \sum_{i=1}^\infty 1_{\left\{\tau_i\le t\right\}}\xi_i$$
and has RCLL paths. Let $\lambda>0$ be a fixed {\it discount factor}.

\begin{definition}\label{d2} A control policy $(T,\Xi)$ is called {\it admissible} if the corresponding discounted controlled cash on hand process is bounded, that is, there exists a constant $K\in(0,\infty)$ such that
$$e^{-\lambda t}|X_t|<K,\quad\hbox{a.s.},\ t\in[0,\infty).$$
\end{definition}

\noindent We emphasize that the constant $K$ is not  "universal"; for different policies it may be different.

 Let the opportunity cost function be $\phi(x)\ge 0$ for the cash on hand  being at level $x$. We shall assume that for some constants $K_1,K_2>0$
\begin {equation}K_1\phi(x)+K_2\ge |x|,\quad x\in\Re.\label{new2}
\end{equation}
 The manager's total cost is the sum of his/her  expected discounted opportunity costs as well as impulse control costs and is given by
\begin{equation}I(x,T,\Xi)=E_x\left[\int_0^\infty e^{-\lambda t}\phi(X_t)dt +\sum_{n=1}^\infty e^{-\lambda\tau_n }g(\xi_n)\right]\label{(A1)}\end{equation}
where the manager's impulse control costs are given by
\begin{equation*}
g(\xi)=
\begin{cases}
C+c\xi, &\textrm{if~} \xi>0, \\
                0, &\textrm{if~} \xi=0, \\
                D-d\xi,&\textrm{if~} \xi<0.
\end{cases}
\end{equation*}
 We assume  that the fixed costs $C$ and $D$ are positive and the variable costs $c$ and $d$ are non-negative constants. The term $e^{-\lambda\tau_n}$ in (\ref{(A1)}) is well defined once we set $e^{-\infty}=0.$

One of our primary objectives in this paper is to identify the optimal impulse control $(T,\Xi)$ that minimizes the above total cost $I(x,T,\Xi)$. The value function of this optimization problem is
$$V(x) = \inf \left\{I(x,T,\Xi),\ (T,\Xi)\ \hbox{is an admissible impulse control}\right\}. $$
Notice  that extending policies according to our convention (\ref{new1}) does not change the corresponding cost.

In the section that follows we show that the impulse control takes the form of a control band policy
which arise frequently as the solution to impulse control problems. Moreover, in Section 5  we provide an
example in which we are able to explicitly identify the parameters corresponding to this policy.

\section{Main Results}
\label{Sec::Main}

In this Section, we provide the main result of the paper, Theorem \ref{t2}, showing that a solution of an ordinary differential equation that satisfies some additional conditions,
must be the value function $V$.  Also included in the statement of Theorem \ref{t2} is the optimal impulse control policy which  turns out
to be a double bandwidth control policy. We begin first with some preliminary results before providing the statement of Theorem \ref{t2}.

For a function $f:\Re\mapsto\Re$ we define the operator
$$Mf(x)=\inf\left\{f(x+\eta)+g(\eta),\ \eta\in\Re\setminus\{0\}\right\}.$$
We shall also use the linear operator ${\cal A}$ associated with the uncontrolled process $Y$, that is, for $f\in C^2(\Re)$
$${\cal A}f(x)={1\over 2}\sigma^2f''(x) + \mu f'(x) +\int_{\Re}\left[f(x+y)-f(x) - f'(x)y1_{\left\{|y|<1\right\}}\right]\nu(dy).$$
Our assumption (\ref{Main}) and Taylor's theorem implies that  the integral on the right-hand side is finite whenever $f'$ and $f''$ are bounded.
Indeed,
\begin{eqnarray}
&&\left|\int_{\Re}\left[f(x+y)-f(x) - f'(x)y1_{\left\{|y|<1\right\}}\right]\nu(dy)\right| \nonumber\\
&\le&\int_{\Re}\left|f(x+y)-f(x)\right|1_{\left\{|y|\ge 1\right\}}\nu(dy)+
\int_{\Re}\left|f(x+y)-f(x)-f'(x)y\right|1_{\left\{|y|< 1\right\}}\nu(dy) \nonumber\\
&\le& \mathrm{const}\times\left[\int_\Re |y|1_{\{|y|\ge 1\}}\nu(dy) + \int_\Re y^21_{\{|y|< 1\}}\nu(dy)\right] \nonumber\\
&<&\infty.\label{new6}
\end{eqnarray}

In order to
prove our results we shall actually need to
extend the domain of ${\cal A}$ to a larger class of functions ${\cal D}$ defined below.

\begin{definition}\label{d1} Let ${\cal D}$ be the class of functions $f:\Re\mapsto\Re$ for which there exist an integer $n\ge 0$ and a set of real numbers $S=\{x_1,x_2,\dots,x_n\}$ (if $n=0$ then $S$ is the empty set) such that the following conditions hold:\hfill\break
(i) $f\in C^1(\Re)\cap C^2(\Re\setminus S)$\hfill\break
(ii) The derivative $f'$ is bounded on $\Re$ and the second derivative $f''$ is bounded on $\Re\setminus S$. We shall call the points in $S$ the exceptional points.
\end{definition}

One can see easily that (\ref{new6}) holds even if $f\in{\cal D}$, thus we
 extend the operator ${\cal A}$ to ${\cal D}$. This way ${\cal A}f(x)$ may be undefined if $x$ is an exceptional point of $f$, but this will not cause any problems. For $f\in {\cal D}$ it may be shown that Ito's rule applied to $f(X_t)$ holds in its usual form (see the Appendix).

We now conjecture that the optimal impulse control policy takes a double bandwidth control policy form. In particular, we assume that there exist constants $a<\alpha\le\beta<b$ such that \begin{equation}\tau^*_n=\inf\left\{t\ge\tau^*_{n-1}:\ X_{t-}+\Delta Y_t\in\Re\setminus(a,b)\right\}\label{(A2)}\end{equation}
(recall that $\tau^*_0=0$)
and that for $n\ge 1$ the jump size is given by
\begin{equation}
\xi^*_n=
\begin{cases}
\beta-\left(X(\tau^*_n-)+\Delta Y(\tau^*_n)\right),&\textrm{if}~X(\tau^*_n-)+\Delta Y(\tau^*_n)\ge b ,\\
           \alpha-\left(X(\tau^*_n-)+\Delta Y(\tau^*_n)\right), &\textrm{if}~X(\tau^*_n-)+\Delta Y(\tau^*_n)\le a.
           \label{(A3)}
\end{cases}
\end{equation}
Note that $\tau_1^*=0$ if and only if $x\in\Re\setminus(a,b)$. For $n>1$ it is possible that $\Delta Y(\tau_n^*)=0$ in which case $X(\tau_n^*-)$ is either equal to $a$ or $b$, $\xi_n^*$ is $\alpha-a$ or $\beta -b$
and $X(\tau_n^*)$ equal to $\alpha$ or  $\beta$, respectively. However, it is also possible that $\Delta Y(\tau_n^*)\not=0$ in which case $X(\tau_n^*-)+\Delta Y(\tau_n^*)$ may be either larger than  $a$ or smaller then $b$, but we still have $X(\tau_n^*)$ equal to $\alpha$ or  $\beta$, respectively.  This control policy is admissible since the corresponding cash at hand process is bounded, even without discounting.

\begin{proposition}\label{p1}
Suppose that for some $a<\alpha\le \beta<b$ the sequence of stopping times $\left(\tau^*_n\right)_{\{n\ge 1\}}$ is  given in (\ref{(A2)}).
If $Y$ is not constant then $\tau^*_n<\infty$, $\tau^*_n<\tau^*_{n+1}$ for $n\ge 1$ almost surely, and $\lim_{n\to\infty}\tau^*_n=\infty$.
\end{proposition}

{\bf Proof:} It is obvious that $\tau^*_n<\infty$ since otherwise the Levy process $Y$ would be bounded, that is, a constant. Next we show
$\tau^*_n<\tau^*_{n+1}$. Let $Y^\alpha$ be the process $Y$ started at $Y_0=\alpha$  and $Y^\beta$ be the process $Y$ started at $Y_0=\beta$.
Let $\tau_{\alpha}=\inf\left\{s\ge 0:\ Y_s^\alpha\notin (a,b)\right\}$ and $\tau_{\beta}=\inf\left\{s\ge 0:\ Y_s^\beta\notin (a,b)\right\}$. By the
right-continuity of $Y$ we have $\tau_{\alpha}>0$ and $\tau_{\beta}>0$ a.s. But conditionally on $X(\tau^*_n)=\alpha$ the inter-arrival time
$\tau^*_{n+1}-\tau^*_n$ has the same distribution as $\tau_\alpha$ and conditionally on $X(\tau^*_n)=\beta$ the inter-arrival time $\tau^*_{n+1}-\tau^*_n$ has the same distribution as $\tau_\beta$. Hence  $\tau^*_n<\tau^*_{n+1}$. Finally we show that  $\lim_{n\to\infty}\tau^*_n=\infty$. We first show that either $P[X(\tau^*_n)=\alpha,\ \hbox{i.o.}]=1$ or $P[X(\tau^*_n)=\beta,\ \hbox{i.o.}]=1$ or both. First note that
if  both $P[X(\tau^*_1)=\alpha | X(0)=\alpha]=1$ and $P[X(\tau^*_1)=\beta | X(0)=\beta]=1$ (which actually is not even possible), then by the strong Markov property of $X_t$, either $P[X(\tau^*_n)=\alpha,\ \hbox{i.o.}]=1$ or $P[X(\tau^*_n)=\beta,\ \hbox{i.o.}]=1$, but not both. On the other hand, if $P[X(\tau^*_1)=\alpha | X(0)=\alpha]<1$ and $P[X(\tau^*_1)=\beta | X(0)=\beta]<1$, then, by the strong Markov property of $X_t$, both $P[X(\tau^*_n)=\alpha,\ \hbox{i.o.}]=1$ and $P[X(\tau^*_1)=\beta,\ \hbox{i.o.}]=1$. Finally, note that if $P[X(\tau^*_1)=\alpha | X(0)=\alpha]=1$ and $P[X(\tau^*_1)=\beta | X(0)=\beta]<1$, then by the strong Markov property, $P[X(\tau^*_n)=\alpha,\ \hbox{i.o.}]=1$. In the reverse case, one has that $P[X(\tau^*_n)=\beta,\ \hbox{i.o.}]=1$
Suppose now, for the sake of argument, that $P[X(\tau^*_n)=\alpha,\ \hbox{i.o.}]=1$. Then let $S_0=0$ and $S_{n+1}=\min\{\tau^*_m>S_n:\ X(\tau^*_m)=\alpha\}$. Then $S_2-S_1,S_3-S_2,\dots$ is an i.i.d. sequence of random variables. Since $S_{n+1}-S_n>0$, a.s., so $ E[S_{n+1}-S_n]>0$ which implies that $\sum _{n=1}^\infty (S_{n+1}-S_n)$ diverges, i.e. $S_n\to\infty$. But this implies $\tau^*_n\to\infty$.

The following is now the main result of this Section.

\begin{theorem}\label{t2} Suppose that there exist constants $a<\alpha\le\beta<b$ and a function $f:\Re\mapsto (0,\infty)$ such that $f\in C^1(\Re)\cap C^2\left(\Re\setminus\{a,b\}\right)$, the second derivative $f''$ is bounded on $\Re\setminus\{a,b\}$ and the following conditions are satisfied:\hfill\break
(i) ${\cal A}f(x) -\lambda f(x) +\phi(x)=0$\quad for $x\in(a,b)$;\hfill\break
(ii) $f(x)\le Mf(x)$\quad for $x\in(a,b)$;\hfill\break
(iii) ${\cal A}f(x) -\lambda f(x) +\phi(x)\ge 0$\quad for $x\in\Re\setminus[a,b]$;\hfill\break
(iv) $f(a)=Mf(a)=f(\alpha)+C+c(\alpha -a),\quad f(b)=Mf(b)=f(\beta)+D+d(b-\beta)$;\hfill\break
(v) $f$ is linear on $(-\infty,a]$ with slope $-c$, and also linear on $[b,\infty)$ with slope $d$.

Then $f(x)=V(x)$. Furthermore, the  control $(T^*,\Xi^*)$ given in {\rm (\ref{(A2)}) - (\ref{(A3)})} with these values for $a,\alpha,\beta, b$ is optimal.
\end{theorem}

In order to provide the proof of Theorem \ref{t2}, we first need the following.

\begin{lemma}\label{l3} Let $(T,\Xi)$  be a control such that $I(x,T,\Xi)<\infty$ Then
\begin{equation}\lim_{t\to\infty}E_x\left[e^{-\lambda t}|X_{t}|\right]=0\label{(1.3)}\end{equation}
and
\begin{equation}\lim_{n\to\infty}\tau_n(\omega)=\infty,\quad a.s.\label{(1.1)}\end{equation}
\end{lemma}

{\bf Proof:} First we prove (\ref{(1.3)}). By (\ref{new2})
$$\int_0^\infty E_x\left[e^{-\lambda t}|X_t|\right]dt \le \int_0^\infty E_x\left[e^{-\lambda t}\left(K_1\phi(X_t)+K_2\right)\right]dt,$$
which is finite by our assumption that $I(x,T,\Xi)<\infty$.
Next we prove (\ref{(1.1)}). Suppose the opposite, i.e., that $P(G)>0$ where $G=\{\omega\in \Omega:\ \tau_n\uparrow\tau<\infty\}$.
Then
\begin{eqnarray*}
&&\sum_{n=1}^\infty E_x\left[e^{-\lambda\tau_n}g(\xi_n)\right] \ge \sum_{n=1}^\infty E_x\left[e^{-\lambda\tau_n}1_G\right]\hbox{min}\left\{C,D\right\}\ge\sum_{n=1}^\infty E_x\left[e^{-\lambda\tau}1_G\right]\hbox{min}\left\{C,D\right\}=\infty
\end{eqnarray*}
which contradicts our assumption that $I(x,T,\Xi)<\infty$.\hfill$\bullet$

By Lemma \ref{l3} without any loss of generality we can and shall consider only those policies for which (\ref{(1.3)}) and (\ref{(1.1)}) hold.



{\bf Proof of Theorem \ref{t2}:} Using condition (v) we can extend (iv) to $(-\infty,a]$ and $[b,\infty)$:
 \begin{equation}f(x)=Mf(x)=f(\alpha)+C+c(\alpha -x),\quad x\le a\label{(1)}\end{equation}
\begin{equation}f(x)=Mf(x)=f(\beta)+D+d(x-\beta),\quad x\ge b.\label{(2)}\end{equation}
Indeed, suppose that $x<a$. Then for $\eta\in[0,a-x]$ the quantity $f(x+\eta)+C+c\eta=f(x)+C$ does not depend on $\eta$. In addition  $f(x+\eta)+g(\eta)$ is a decreasing function of $\eta$ on $(-\infty,0)$. Hence by (iv) and (v)
$$Mf(x)=\inf\{f(x+\eta)+C+c\eta,\ \eta\ge a-x\}=\inf\{f(a+\gamma)+C+c(a+\gamma-x);\gamma\ge 0\}$$
$$=Mf(a)+c(a-x)=f(a)+c(a-x)=f(x)$$
Also by (iv)
$$Mf(a)+c(a-x)=f(\alpha)+C+c(\alpha-a)+c(a-x)=f(\alpha)+C+c(\alpha-x).$$
We deal with the $x>b$ case similarly. Let $(T,\Xi)=\{\tau_1,\tau_2,\dots,\xi_1,\xi_2,\dots\}$ be an arbitrary impulse control for which (\ref{(1.3)}) and (\ref{(1.1)}) hold.
We define the local martingale
\begin{equation}U_t=\int_{(0,t]\times\Re}   e^{-\lambda s}\left\{f(X_{s-}+y)-f(X_{s-})-f'(X_{s-})y1_{\{|y|<1\}}\right\}\left(J(ds,dy)-ds\nu(dy)\right)\label{new8}\end{equation}
Since $U$ is a local martingale, there exists a sequence of stopping times $\left(S_m\right)$ such that $\lim_{m\to\infty}S_m=\infty$, almost surely, and $\{U_{t\wedge S_m},\ t\in[0,\infty)\}$ is a martingale for every $m$.
Fix $t>0$ $m\ge 1$. For the sake of brevity let $T_n=\min\{\tau_n,S_m,t\}$. The stopping time $T_n$ depends on $m$ and $t$ but we suppress these dependencies in the notation.
We have the following decomposition
\begin{eqnarray}
&&e^{-\lambda T_n}f\left(X_{T_n}\right)-f(x)  \label{(3)}\\
&=&\sum_{i=1}^n\left\{e^{-\lambda T_i}f\left(X_{T_i-}+\Delta Y_{T_i}\right) - e^{-\lambda T_{i-1}}f\left(X_{T_{i-1}}\right)\right\}+\sum_{i=1}^ne^{-\lambda T_i}\left[f\left(X_{T_i}\right)-f\left(X_{T_i-}+\Delta Y_{T_i}\right)\right]. \nonumber
\end{eqnarray}
First we deal with the first sum on the right-hand side of (\ref{(3)}).
By the generalized Ito's rule (Proposition \ref{p6})  and integration by parts   follows that
\begin{eqnarray*}
&&e^{-\lambda T_i}f\left(X_{T_i-}+\Delta Y_{T_i}\right) - e^{-\lambda T_{i-1}}f\left(X_{T_{i-1}}\right) \\
&=&\int_{(T_{i-1},T_i]}e^{-\lambda s} f'(X_{s-})\left\{dA_s+dM_s+\sigma dw_s\right\}\\
&&+\int_{T_{i-1}}^{T_i}e^{-\lambda s}\left\{{\sigma^2\over 2}f''(X_{s-})+\mu f'(X_{s-})-\lambda f(X_{s-})\right\}ds \\
&&+\sum_{T_{i-1}<s\le T_i}e^{-\lambda s}\left\{f(X_s)-f(X_{s-})-f'(X_{s-})\Delta X_s\right\}.
\end{eqnarray*}
The right-hand side can be written as
\begin{eqnarray*}
&&\int_{(T_{i-1},T_i]}e^{-\lambda s} f'(X_{s-})\left\{dM_s+\sigma dw_s\right\} \\
&&+\int_{T_{i-1}}^{T_i}e^{-\lambda s}\left\{{\sigma^2\over 2}f''(X_{s-})+\mu f'(X_{s-})-\lambda f(X_{s-})\right\}ds\\
&&+\sum_{T_{i-1}<s\le T_i}e^{-\lambda s}\left\{f(X_s)-f(X_{s-})-f'(X_{s-})\Delta X_s1_{|\Delta X_s|<1}\right\},
\end{eqnarray*}
which is equal to
\begin{eqnarray*}
&&\int_{(T_{i-1},T_i]}e^{-\lambda s} f'(X_{s-})\left\{dM_s+\sigma dw_s\right\} \\
&&+\int_{T_{i-1}}^{T_i}e^{-\lambda s}\left\{{\sigma^2\over 2}f''(X_{s-})+\mu f'(X_{s-})-\lambda f(X_{s-})\right\}ds+ \\
&&+\int_{(T_{i-1},T_i]\times\Re}   e^{-\lambda s}\left\{f(X_{s-}+y)-f(X_{s-})-f'(X_{s-})y1_{\{|y|<1\}}\right\}J(ds,dy).
\end{eqnarray*}
Since the left-hand side of (\ref{new6}) is finite,
the above expression can be written as
\begin{eqnarray*}
&&\int_{(T_{i-1},T_i]}e^{-\lambda s} f'(X_{s-})\left\{dM_s+\sigma dw_s\right\}\\
&&+\int_{T_{i-1}}^{T_i}e^{-\lambda s}\left\{{\sigma^2\over 2}f''(X_{s-})+\mu f'(X_{s-})-\lambda f(X_{s-})\right\}ds\\
&&+\int_{(T_{i-1},T_i]\times\Re}   e^{-\lambda s}\left\{f(X_{s-}+y)-f(X_{s-})-f'(X_{s-})y1_{\{|y|<1\}}\right\}\left(J(ds,dy)-ds\nu(dy)\right)\\
&&+\int_{(T_{i-1},T_i]\times\Re}   e^{-\lambda s}\left\{f(X_{s-}+y)-f(X_{s-})-f'(X_{s-})y1_{\{|y|<1\}}\right\}ds\nu(dy)\\
&=&\int_{(T_{i-1},T_i]}e^{-\lambda s} f'(X_{s-})\left\{dM_s+\sigma dw_s\right\} \\
&&+\int_{T_{i-1}}^{T_i}e^{-\lambda s}\left\{{\cal A}(X_{s-})-\lambda f(X_{s-})\right\}ds\\
&&+\int_{(T_{i-1},T_i]\times\Re}   e^{-\lambda s}\left\{f(X_{s-}+y)-f(X_{s-})-f'(X_{s-})y1_{\{|y|<1\}}\right\}\left(J(ds,dy)-ds\nu(dy)\right).
\end{eqnarray*}
By conditions (i) and (iii) we end up with the inequality
\begin{eqnarray}
&&e^{-\lambda T_i}f\left(X_{T_i-}+\Delta Y_{T_i}\right) - e^{-\lambda T_{i-1}}f\left(X_{T_{i-1}}\right)  \nonumber \\
&\ge&\int_{(T_{i-1},T_i]}e^{-\lambda s} f'(X_{s-})\left\{dM_s+\sigma dw_s\right\}  \nonumber\\
&&+\int_{(T_{i-1},T_i]\times\Re}   e^{-\lambda s}\left\{f(X_{s-}+y)-f(X_{s-})-f'(X_{s-})y1_{\{|y|<1\}}\right\}\left(J(ds,dy)-ds\nu(dy)\right)  \nonumber\\
&&-\int_{T_{i-1}}^{T_i}e^{-\lambda s}\phi(X_s)ds.\label{(3.1)}
\end{eqnarray}
There is a minor problem here since  (i) and (iii) implies the inequality ${\cal A}f(x) -\lambda f(x) +\phi(x)\ge 0$ only for $x\in\Re\setminus\{a,b\}$.
But either $\sigma=0$, in which case the inequality holds for all $x\in\Re$, or, $\sigma\not=0$, in which case (\ref{(3.1)}) still holds by Lemma \ref{l8} in the Appendix.

For our candidate optimal policy $(T^*,\Xi^*)$ actually equality holds in (\ref{(3.1)}) by condition (i).

Next we deal with the second sum on the right-hand side of (\ref{(3)}). By  (\ref{(1)}), (\ref{(2)}), and (ii) whenever $T_i=\tau_j$ for some $j\ge 1$ then
\begin{equation}f(X_{T_i}) - f\left(X_{T_i-}+\Delta Y_{T_i}\right)\ge -g\left(X_{T_i}-(X_{T_i-}+\Delta Y_{T_i})\right) = -g(\xi_j).\label{new4}\end{equation}
 If $T_n=t\wedge S_m$ then the left-hand side of (\ref{new4}) is zero, hence
\begin{equation}\sum_{i=1}^ne^{-\lambda T_i}\left[f\left(X_{T_i}\right)-f\left(X_{T_i-}+\Delta Y_{T_i}\right)\right]\ge -\sum_{j\le n,\tau_j\le t\wedge S_m} e^{-\lambda\tau_j}g(\xi_j)\label{new5}\end{equation}
By (\ref{(1)}) and (\ref{(2)}) this becomes an  equality for $(T^*,\Xi^*)$.

Adding up (\ref{(3.1)}) and (\ref{new5}) and considering (\ref{new8}), (\ref{(3)}), we get
\begin{eqnarray*}
&&e^{-\lambda T_n}f\left(X_{T_n}\right)-f(x)\\
&\ge&\int_{(0,T_n]}e^{-\lambda s} f'(X_{s-})\left\{dM_s+\sigma dw_s\right\}+U_{T_n}
-\int_{0}^{T_n}e^{-\lambda s}\phi(X_s)ds-
\sum_{j\le n,\tau_j\le t\wedge S_m}e^{-\lambda\tau_j}g(\xi_j).
\end{eqnarray*}
After taking expectations on both sides, all martingale terms will disappear. Indeed, the local martingale
$$Z_t=\int_{(0,t]}e^{-\lambda s}f'(X_{s-})dM_s$$
has quadratic variation
$$\left[Z\right]_t=\int_{(0,t]} e^{-2\lambda s}  \left(f'(s)\right)^2d\left[M\right]_s$$
(\cite{protter}, Chapter II, Theorem 29), and the boundedness of $f'$ and (\ref{new7}) imply that $Z$ is a square-integrable martingale (see also \cite{protter}, Chapter II, Corollary 3 to Theorem 27). Then by the Optional sampling Theorem for bounded stopping times we have $E[Z_{T_n}]=0$.
Also, since $\left\{U_{s\wedge S_m},s<\infty\right\}$ is a martingale, another application of the Optional sampling Theorem implies that $E[U_{\tau_n\wedge t\wedge S_m}]=E[U_{T_n}]=0$.
 Hence after taking expectations on both sides  we end up with
$$-E_x\left[e^{-\lambda T_n}f\left(X_{T_n}\right)\right]+f(x)\le
E_x\left[\int_{0}^{T_n}e^{-\lambda s}\phi(X_s)ds+\sum_{j:\tau_j\le T_n}e^{-\lambda\tau_j }g(\xi_j)\right].$$
Now let $n\to \infty$. By the continuity of $f$ and condition (v) $f$ satisfies a linear growth condition $f(x)\le K_3|x|+K_4$ for some positive constants $K_3,K_4$, hence
$0\le  e^{-\lambda T_n}f\left(X_{T_n}\right)\le e^{-\lambda T_n}(K_3\left|X_{T_n}\right|+K_4)$, which is bounded by the admissibility of $(T,\Xi)$. We apply the Dominated Convergence Theorem on the left-hand side and the Monotone Convergence Theorem on the right-hand side, and get the inequality
$$-E_x\left[e^{-\lambda (t\wedge S_m)}f\left(X_{t\wedge S_m}\right)\right]+f(x)\le
E_x\left[\int_{0}^{t\wedge S_m}e^{-\lambda s}\phi(X_s)ds+\sum_{j:\tau_j\le t\wedge S_m}e^{-\lambda\tau_j}g(\xi_j)\right].$$
We apply exactly the same method as $m\to\infty$ and get
$$-E_x\left[e^{-\lambda t}f\left(X_t\right)\right]+f(x)\le
E_x\left[\int_{0}^{t}e^{-\lambda s}\phi(X_s)ds+\sum_{j:\tau_j\le t}e^{-\lambda\tau_j}g(\xi_j)\right].$$
The linear growth of $f$ and (\ref{(1.3)}) imply that $E\left[e^{-\lambda t}f\left(X_t\right)\right]\to\infty$ as $t\to \infty$. So by letting $t\to\infty$ another application of the Monotone Convergence Theorem on the right-hand side gives
$$f(x)\le E_x\left[\int_{0}^{\infty}e^{-\lambda s}\phi(X_s)ds+\sum_{j=1}^\infty e^{-\lambda\tau_j }g(\xi_j)\right]$$
Equality holds $(T^*,\Xi^*)$ hence the proof is complete.\hfill$\bullet$

\section{Analysis of the Optimal Control $(T^*,\Xi^*)$}
\label{Sec::Analysis}

We now set out to determine the transient and steady-state behavior of the controlled cash on hand process $X_t$
under the optimal control $(T^*,\Xi^*)$ assuming that $Y_t$ is a spectrally positive Levy process. In other words,
\begin{eqnarray*}
\int_{\Re_{-}} \nu(dy)=0
\end{eqnarray*}
and $Y_t$ is not a subordinator. The case of a spectrally negative Levy process may be treated similarly. Our main result will be to determine the Laplace transform (with respect to time) of the transition probabilities of $X_t$ and also to determine the limiting
distribution of $X_t$.

We begin by determining the Laplace transform of  the transition probabilities of $X_t$.  Let $\mathcal{A} \in \mathcal{B}(\Re)$ be a Borel set of $\Re$ and let $e_q$ be an exponential random variable with
rate $q$ independent of $X$. Now consider $  P_x[X_{e_q} \in \mathcal{A}]$ for $x \in (a,b)$. It then follows conditioning on the value of $e_q$ relative to the stopping time $\tau^*_1$ that
\begin{equation}P_x[X_{e_q} \in {\cal A}] = E_x[ 1 \{X_{e_q} \in {\cal A}\} 1\{e_q < \tau^*_1 \}]+E_x[1\{X_{e_q} \in {\cal A}\} 1\{e_q \geq \tau^*_1 \}]. \label{J18}
\end{equation}
However, since $X_t=Y_t$ for $0 \leq t < \tau^{*}_1$, it follows by the memoryless property of the exponential distribution and the strong Markov property that $E_x[1\{X_{e_q} \in {\cal A}\}1\{e_q \geq \tau^*_1\}] $ is equal to
\begin{eqnarray*}
E_{\alpha}[1\{X_{e_q} \in {\cal A}\}]P_x[\{e_q \geq \tau^*_1 \} \cap \{ Y(\tau^*_1) \leq a\}]+E_{\beta}[1\{X_{e_q} \in {\cal A}\}] P_x[\{e_q \geq \tau^*_1\} \cap \{    Y(\tau^*_1) \geq b\}].
\end{eqnarray*}
Substituting the above into (\ref{J18}), we have
\begin{eqnarray}
E_x[1\{X_{e_q} \in {\cal A}\}]&=&E_x[1\{X_{e_q} \in {\cal A}\}1\{e_q < \tau^*_1\}] \label{Pre}\\
&&+ E_{\alpha}[1\{X_{e_q} \in {\cal A}\}]P_x[\{e_q \geq \tau^*_1\} \cap \{ Y(\tau^*_1) \leq a\}]  \nonumber \\
&&+E_{\beta}[1\{X_{e_q} \in {\cal A}\}] P_x[\{e_q \geq \tau^*_1\} \cap \{ Y(\tau^*_1) \geq b\}]. \nonumber
\end{eqnarray}
It therefore remains to determine expressions for the three quantities on the righthand side of (\ref{Pre}).
We proceed term by term.

First note that $\tau^{*}_1$ is equal to the fist time the Levy process $Y_t$ exits the open interval $(a,b)$. We then have that in general one may write
\begin{eqnarray}
E_x[1\{X_{e_q} \in {\cal A}\}1\{e_q < \tau^*_1\}]&=&E_x[1\{Y_{e_q} \in {\cal A}\}1\{e_q < \tau^*_1\}] \label{Ex}\\
&=&q\int_0^{\infty} e^{-qt} E_{x}[1\{Y_{t} \in {\cal A}\}1\{ \tau^*_1 > t\}]dt \nonumber\\
&=&q \int_{{\cal A}} \int_0^{\infty} e^{-qt} P_{x}[Y(t) \in dy,  \tau^*_1 > t]dt \nonumber \\
&=&q \int_{{\cal A}} U^{(q)}(x,dy), \nonumber
\end{eqnarray}
where $U^{{(q)}}$ is the $q$-potential measure of $Y_t$. This then provides an expression for the first term on the righthand side of (\ref{Pre}). Also note, by Theorem 8.7 of \cite{K06}, if $Y_t$ is spectrally positive then its $q$-potential measure $U^{(q)}(x,dy)$ has a density $u^{(q)}(x,y)$ given by
\begin{eqnarray}
u^{(q)}(x,y)&=&W^{(q)}(b-x)\frac{W^{(q)}(y-a)}{W^{(q)}(b-a)} - W^{(q)}(y-x), \label{density::U}
\end{eqnarray}
where
\begin{eqnarray}
\int_0^{\infty}e^{-sy}W^{(q)}(y)dy &=& \frac{1}{\psi(-s)-q}, \label{LaplaceW}
\end{eqnarray}
whenever $-s$ is large enough so that $\psi(-s) > q$ and $\psi(s)=\log E_0[e^{sY_1}]$ is the Laplace exponent of $Y_t$. Note that $\psi(s) < \infty$ for all $s \leq 0$ by the spectral positivity of $Y$.

Next note that setting $x=\alpha$ in (\ref{Pre}), we obtain
\begin{eqnarray}
&&E_{\alpha}[1\{X_{e_q} \in {\cal A}\}](1 -P_{\alpha}[\{e_q \geq \tau^*_1\} \cap \{ Y(\tau^*_1) \leq a\}) ) \label{SimilarOne}\\
&&-E_{\beta}[1\{X_{e_q} \in {\cal A}\}] P_{\alpha}[\{e_q \geq \tau^*_1\} \cap \{ Y(\tau^*_1) \geq b\}] \nonumber \\
&=&E_{\alpha}[1\{X_{e_q} \in {\cal A}\}1\{e_q < \tau^*_1\}],  \nonumber
\end{eqnarray}
and similarly, setting $x=\beta$, we have
\begin{eqnarray}
&&E_{\beta}[1\{X_{e_q} \in {\cal A}\}](1 -P_{\beta}[\{e_q \geq \tau^*_1\} \cap \{Y(\tau^*_1) \geq b\}] ) \label{SimilarTwo}\\
&&-E_{\alpha}[1\{X_{e_q} \in {\cal A}\}] P_{\beta}[\{e_q \geq \tau^*_1\} \cap \{ Y(\tau^*_1) \leq a\}] \nonumber \\
&=&E_{\beta}[1\{X_{e_q} \in {\cal A}\}1\{e_q < \tau^*_1\}].  \nonumber
\end{eqnarray}
(\ref{SimilarOne}) and (\ref{SimilarTwo}) constitute a set of linear equations for
$E_{\alpha}[1\{X_{e_q} \in {\cal A}\}]$ and $E_{\beta}[1\{X_{e_q} \in {\cal A}\}]$.
Moreover,   so long as $q > 0$, we have that
\begin{eqnarray}
&&(1 -P_{\beta}[\{e_q \geq \tau^*_1\} \cap \{ Y(\tau^*_1) \geq b\}] )(1 -P_{\alpha}[\{e_q \geq \tau^*_1\} \cap \{Y_{\tau^*_1} \leq a\}] ) \nonumber \\
&>& P_{\alpha}[\{e_q \geq \tau^*_1\} \cap \{Y(\tau^*_1) \geq b\}]P_{\beta}[\{e_q \geq \tau^*_1\} \cap \{Y_{\tau^*_1} \leq a\}]  \nonumber
\end{eqnarray}
and so the determinant associated with (\ref{SimilarOne}) and (\ref{SimilarTwo}) is non-zero and hence a
solution exists. Solving for $E_{\alpha}[1\{X_{e_q} \in {\cal A}\}]$ and $E_{\beta}[1\{X_{e_q} \in {\cal A}\}]$
then yields that $E_{\alpha}[1\{X_{e_q} \in {\cal A}\}]$ is given by
\begin{eqnarray}
{1\over C_{a,\alpha,\beta,b}}\times
\left( \begin{array}{c} E_{\alpha}[1\{X_{e_q} \in {\cal A}\}1\{e_q < \tau^*_1\}](1 -P_{\beta}[\{e_q \geq \tau^*_1\} \cap \{ Y(\tau^*_1) \geq b\}]    \\\
 +E_{\beta}[1\{X_{e_q} \in {\cal A}\}1\{e_q < \tau^*_1\}]P_{\alpha}[\{e_q \geq \tau^*_1\} \cap \{ Y(\tau^*_1) \geq b\}]      \label{J28}
\end{array}
\right)
\end{eqnarray}
and $E_{\beta}[1\{X_{e_q} \in {\cal A}\}]$ is given by
\begin{eqnarray}
{1\over C_{a,\alpha,\beta,b}}\times
\left(
 \begin{array}{c}
 E_{\alpha}[1\{X_{e_q} \in {\cal A}\}1\{e_q < \tau^*_1\}]P_{\beta}[\{e_q \geq \tau^*_1\} \cap \{Y_{\tau^*_1} \leq a\}]     \\
 +E_{\beta}[1\{X_{e_q} \in {\cal A}\}1\{e_q < \tau^*_1\}](1-P_{\alpha}[\{e_q \geq \tau^*_1\} \cap \{Y_{\tau^*_1} \leq a\}])
\end{array}
\right),
 \label{J29}
\end{eqnarray}
where
\begin{eqnarray}
C_{a,\alpha,\beta,b}&=&(1 -P_{\beta}[\{e_q \geq \tau^*_1\} \cap \{ Y(\tau^*_1) \geq b\}] )(1 -P_{\alpha}[\{e_q \geq \tau^*_1\} \cap \{Y_{\tau^*_1} \leq a\}] )  \nonumber\\
&&- P_{\alpha}[\{e_q \geq \tau^*_1\} \cap \{Y(\tau^*_1) \geq b\}]P_{\beta}[\{e_q \geq \tau^*_1\} \cap \{Y_{\tau^*_1} \leq a\}]. \nonumber
\end{eqnarray}
We now proceed to compute the terms appearing on the right-hand sides of (\ref{J28}) and (\ref{J29}).

Firs note that the expressions of the form $E_{x}[1\{X_{e_q} \in {\cal A}\}1\{e_q < \tau^*_1\}]$  appearing in
(\ref{J28}) and (\ref{J29}) have already been determined by (\ref{Ex}). It remains to determine expressions for
the probabilities appearing in (\ref{J28}) and (\ref{J29}). However, note that
\begin{eqnarray*}
P_x[\{e_q \geq \tau^*_1\} \cap \{Y_{\tau^*_1} \leq a\}]&=& \int_0^{\infty}qe^{-qt}P_x[\{t \geq \tau^*_1\} \cap \{Y_{\tau^*_1} \leq a\}] dt\\
&=&E_x \left[ \int_{\tau^*_1}^{\infty}qe^{-qt}dt 1\{Y_{\tau^*_1} \leq a\}\right] \\
&=&E_x \left[ e^{-q \tau^*_1} 1\{Y_{\tau^*_1} \leq a\} \right] \\
&=& Z^{(q)}(b-x) - Z^{(q)}(b-a)\frac{W^{(q)}(b-x)}{W^{(q)}(b-a) },
\end{eqnarray*}
where the final equality follows from Theorem 8.1 of \cite{K06} and we have the relationship
\begin{eqnarray*}
Z^{(q)}(x) &=& 1 + q \int_0^x W^{(q)}(y)dy.
\end{eqnarray*}
In a similar fashion, using Theorem 8.1 of \cite{K06}, one may compute
\begin{eqnarray*}
P_x[\{e_q \geq \tau^*_1\} \cap \{Y(\tau^*_1) \geq b\}] &=& \frac{W^{(q)}(b-x)}{W^{(q)}(b-a) }.
\end{eqnarray*}
Substituting the above results into (\ref{J28}) and (\ref{J29}), and subsequently into (\ref{Pre}), one obtains an expression for the Laplace transform of the transition probabilities of $X_t$.

We now proceed towards obtaining an expression for the limiting distribution of $X_t$ as $t \rightarrow \infty$.
Our main result in this regard will be the following. Recall the definition of $U^{(q)}$ as the $q$-potential measure associated with $Y$.

\begin{proposition}\label{prop:stat}
If $Y_t$ is spectrally positive, then under a double bandwidth control policy $(a,\alpha,\beta,b)$, the limiting distribution of $X_t$
is given by
\begin{eqnarray*}
\pi(\mathcal{A})&=&\frac{1}{K_{a,\alpha,\beta,b} }\left(\left(1-\frac{W^{(0)}(b-\beta)}{W^{(0)}(b-a)  }\right) \int_{\mathcal{A}}U^{(0)}(\alpha,dy)
 + \left(\frac{W^{(0)}(b-\alpha)}{W^{(0)}(b-a)}    \right) \int_{\mathcal{A}}U^{(0)}(\beta,dy) \right),
\end{eqnarray*}
for each $\mathcal{A} \in \mathcal{B}$, where $K_{a,\alpha,\beta,b}$ is an appropriate normalizing constant, as given by (\ref{Def::K}).
\end{proposition}

We prove Proposition \ref{prop:stat} in a series of lemmas. First note that by the strong Markov property, $X_t$ is a regenerative process with possible regeneration points either $\alpha$ or $\beta$. Let us consider the point $\alpha$ and define $n^{\star}_{\alpha}=\inf \{n \geq 1 : X(\tau^{\star}_n)=\alpha\}$. By the standard theory of regenerative processes, see for instance Theorem 1.2 of Chapter VI of \cite{AS03}, if we may show that $E_{\alpha}[\tau_{n ^{*}_\alpha}] < \infty$ and that $\tau_{n^{*}_\alpha}$
is nonlattice, then  $\lim_{t\to\infty}P_x(X_t\in\cal{A})=\pi({\cal A})$ exists for all ${\cal A}\in{\cal B}(\Re)$ and $x \in \Re$ and is given by
\begin{eqnarray*}
\pi({\cal A}) &=& \frac{E_{\alpha}[\int_0^{\tau_{n^{*}_\alpha}} 1\{X_s\in {\cal A} \} ds]}{E_{\alpha}[\tau_{n^{*}_\alpha} ]}.
\end{eqnarray*}
The following lemma now shows that $E_{\alpha}[\tau_{n^{*}_\alpha}] < \infty$.

\begin{lemma}\label{Prop::Steady}
If the Levy process $Y_t$ is spectrally positive, then
$E_{\alpha}[\tau_{n^{*}_\alpha}] < \infty$.
\end{lemma}

{\bf Proof:} Note first that
\begin{eqnarray*}
\tau_{n^{\star}_{\alpha}} &=& \tau^{\star}_{1} 1\{Y(\tau^{*}_1) \leq a \} +  \tau_{n^{\star}_{\alpha}}  1\{Y(\tau^{*}_1) \geq b \} \\
&=& \tau^{\star}_{1} +  (\tau_{n^{\star}_{\alpha}}- \tau^{*}_1)  1\{Y(\tau^{*}_1) \geq b \} .
\end{eqnarray*}
Hence, by the strong Markov property,
\begin{eqnarray*}
E_{\alpha}[\tau_{n^{\star}_{\alpha}}] &=& E_{\alpha}[\tau^{\star}_{1}] +  E_{\alpha}[(\tau_{n^{\star}_{\alpha}}-\tau^{*}_1)  1\{Y(\tau^{*}_1) \geq b \} ] \\
 &=& E_{\alpha}[\tau^{\star}_{1}] +  E_{\beta}[\tau_{n^{\star}_{\alpha}}]P_{\alpha} [Y(\tau^{*}_1) \geq b ].
\end{eqnarray*}
Similarly, we may show
\begin{eqnarray*}
E_{\beta}[\tau_{n^{\star}_{\alpha}}] &=& E_{\beta}[\tau^{\star}_{1}] +  E_{\beta}[\tau_{n^{\star}_{\alpha}}]P_{\beta} [Y(\tau^{*}_1) \geq b ],
\end{eqnarray*}
from which we obtain
\begin{eqnarray*}
E_{\beta}[\tau_{n^{\star}_{\alpha}}] &=& \frac{E_{\beta}[\tau^{\star}_{1}]}{1-P_{\beta} [Y(\tau^{*}_1) \geq b ] }.
\end{eqnarray*}
Now note that since $Y_t$ is spectrally positive, we have by Theorem 8.1 of \cite{K06} that
\begin{eqnarray*}
P_{\beta} [Y(\tau^{*}_1) \geq b ] &=& \frac{W^{0}(b-\beta)}{W^{0}(b-a)}~<~ 1
\end{eqnarray*}
and so it suffices from the above to show  $E_{\alpha}[\tau^{\star}_{1}],E_{\beta}[\tau^{\star}_{1}] < \infty $.

We now show that in general for $x \in (a,b)$, $E_{x}[\tau^{\star}_{1}] <\infty$.
Recall by \cite{K06}, the potential measure of $Y_t$ upon exiting $[a,b]$ is given by
\begin{eqnarray*}
U(x,dy) = \int_0^{\infty}P_x[Y_t \in dy, \tau^{\star}_1 > t]dt.
\end{eqnarray*}
Integrating over $[a,b]$, we obtain that
\begin{eqnarray*}
\int_{[a,b]}U(x,dy) &=& \int_{[a,b]}\int_0^{\infty}P_x[Y_t \in dy, \tau^{\star}_1 > t]dt \\
&=& \int_0^{\infty} \int_{[a,b]} P_x[Y_t \in dy, \tau^{\star}_1 > t]dt \\
&=& \int_0^{\infty} P_x[ \tau^{\star}_1 > t]dt \\
&=&E_{x}[\tau^{\star}_{1}].
\end{eqnarray*}
However, by Theorem 8.7 of \cite{K06}, since $Y_t$ is spectrally positive, $U(x,dy)$ has a density given by
\begin{eqnarray*}
u(x,y)&=&W(b-x)\frac{W(y-a)}{W(b-a)} - W(y-x).
\end{eqnarray*}
Integrating over $[a,b]$, we therefore find that
\begin{eqnarray*}
\int_{[a,b]}U(x,dy) &=& \int_{[a,b]}\left(W(b-x)\frac{W(y-a)}{W(b-a)} - W(y-x)\right)dy\\
&<& \infty,
\end{eqnarray*}
where the inequality follows since $W$ is bounded on compact sets. By the above, this completes the proof. $\bullet$

The following lemma allows us to take the limit as $q \rightarrow 0$ in (\ref{J28}) and (\ref{J29}) in order to obtain the limiting distribution of $X_t$.

\begin{lemma}\label{pj1} For each $x \in \Re$,
\begin{equation} \lim_{q\to 0}P_x[X(e_q)\in\cal {A}]=\pi({\cal A})\label{j1}\end{equation}
\end{lemma}
{\bf Proof:} Select $T>0$ large enough so that $\left|P_x[X_t\in\cal {A}] - \pi({\cal A})\right|<\epsilon$. Then
$$\left|P_x[X(e_q)\in\cal {A}]-\pi({\cal A})\right|=\left|\int_0^\infty \left\{P_x[X_t\in\cal {A}] -\pi({\cal A})\right\}qe^{-qt}dt\right|$$
$$\le \left|\int_0^T \left\{P_x[X_t\in\cal {A}] -\pi({\cal A})\right\}qe^{-qt}dt\right| + \left|\int_T^\infty \left\{P_x[X_t\in\cal {A}] -\pi({\cal A})\right\}qe^{-qt}dt\right|.$$
The second term in the above expression is bounded by $\epsilon$ and the first term converges to zero as $q\to 0$ by the Dominated Convergence Theorem, which completes the proof. $\bullet$

We now provide the proof of Proposition \ref{prop:stat}.

\noindent {\bf Proof of Proposition \ref{prop:stat}:}
Using Lemma \ref{pj1}, we now wish to take limits $q \rightarrow 0$ in (\ref{J28}) in order to determine
the limiting distribution $\delta$. However, both the numerator and
denominator in (\ref{J28}) converge to 0 as $q \rightarrow 0$ and so we must apply L'Hoptial's rule. Before doing so, however, we first must verify that both the numerator and denominator in (\ref{J28}) are differentiable.

By Lemma 8.3 and Corollary 8.5 in \cite{K06} we have that  for each $x > 0$, both $W^{q}(x)$ and $Z^{(q)}(x)$ are differentiable in $q$. Moreover, since
\begin{eqnarray*}
\frac{W^{(q)}(b-x)}{W^{(q)}(b-a) } &=& E_x[e^{-q \tau^{*}_1}1\{Y_{\tau^{*}_1} \geq b\} ],
\end{eqnarray*}
for $a \leq  x < b$, it follows that
\begin{eqnarray*}
\frac{d}{dq}  \frac{W^{(q)}(b-x)}{W^{(q)}(b-a) }  &=&- E_x[\tau^{*}_1 e^{-q \tau^{*}_1}1\{Y_{\tau^{*}_1} \geq b\} ]\\
                                                  &<&\infty,
\end{eqnarray*}
where the inequality follows as in the proof of Lemma \ref{Prop::Steady}. Finally, since for each $a \leq x \leq b$, $U^{q}(x)$
has a density $u^{(q)}(x,y)$  given by (\ref{density::U}) it follows that for each $\mathcal{A} \in \mathcal{B}(\Re)$,
\begin{eqnarray*}
\frac{d}{dq}\int_{\mathcal{A}}U^{(q)}(x,dy) &=& \int_{\mathcal{A}}\frac{d}{dq} u^{(q)}(x,y)dy.
\end{eqnarray*}
Thus, noting that
\begin{eqnarray}
&&(E_{\alpha}[1\{X_{e_q} \in {\cal A}\}1\{e_q < \tau^*_1\}](1 -P_{\beta}[\{e_q \geq \tau^*_1\} \cap \{ Y(\tau^*_1) \geq b\}] ) \label{FirstWrittenOut} \\
&&+E_{\beta}[1\{X_{e_q} \in {\cal A}\}1\{e_q < \tau^*_1\}]P_{\alpha}[\{e_q \geq \tau^*_1\} \cap \{ Y(\tau^*_1) \geq b\}]) \nonumber \\
&=& q\int_{\mathcal{A}}U^{(q)}(\alpha,dy)\left(1-\frac{W^{(q)}(b-\beta)}{W^{(q)}(b-a)  }\right)
 + q\int_{\mathcal{A}}U^{(q)}(\beta,dy)\left(\frac{W^{(q)}(b-\alpha)}{W^{(q)}(b-a)}    \right)  \nonumber
\end{eqnarray}
and
\begin{eqnarray}
&&((1 -P_{\beta}[\{e_q \geq \tau^*_1\} \cap \{ Y(\tau^*_1) \geq b\}] )(1 -P_{\alpha}[\{e_q \geq \tau^*_1\} \cap \{Y_{\tau^*_1} \leq a\}] )   \label{SecondWrittenOut} \\
&&- P_{\alpha}[\{e_q \geq \tau^*_1\} \cap \{ Y(\tau^*_1) \geq b\}]P_{\beta}[\{e_q \geq \tau^*_1\} \cap \{Y_{\tau^*_1} \leq a\}] ) \nonumber \\
&=&  \left(1-\frac{W^{(q)}(b-\beta)}{W^{(q)}(b-a)  }\right)\left(1-\left(Z^{(q)}(b-\alpha)-Z^{(q)}(b-a)\frac{W^{q}(b-\alpha)}{
W^{q}(b-a)}  \right)   \right)  \nonumber\\
&& -\left(Z^{(q)}(b-\beta)-Z^{(q)}(b-a)\frac{W^{q}(b-\beta)}{
W^{q}(b-a)}  \right)\frac{W^{(q)}(b-\alpha)}{W^{(q)}(b-a)  }   , \nonumber
\end{eqnarray}
we see that both the numerator and denominator in (\ref{J28}) are differentiable.

Let us now take derivatives on the righthand sides of (\ref{FirstWrittenOut}) and (\ref{SecondWrittenOut}).

Taking the derivative of the right hand side of (\ref{FirstWrittenOut}) and evaluating at $q=0$ we obtain
\begin{eqnarray*}
\int_{\mathcal{A}}U^{(0)}(\alpha,dy)\left(1-\frac{W^{(0)}(b-\beta)}{W^{(0)}(b-a)  }\right)
 + \int_{\mathcal{A}}U^{(0)}(\beta,dy)\left(\frac{W^{(0)}(b-\alpha)}{W^{(0)}(b-a)}    \right).
\end{eqnarray*}

Next, recalling that
\begin{eqnarray*}
Z^{q}(x) &=& 1+ q \int_0^x W^{(q)}(y)dy,
\end{eqnarray*}
it follows upon taking the derivative of the righthand side of (\ref{SecondWrittenOut}) and evaluating at $q=0$
that we obtain
\begin{eqnarray}
&&\frac{d}{dq} \left( \left(1-\frac{W^{(q)}(b-\beta)}{W^{(q)}(b-a)  }\right)\left(1-\left(Z^{(q)}(b-\alpha)-Z^{(q)}(b-a)\frac{W^{q}(b-\alpha)}{
W^{q}(b-a)}  \right)   \right) \right)   \nonumber\\
&&\left. -\left(\left(Z^{(q)}(b-\beta)-Z^{(q)}(b-a)\frac{W^{q}(b-\beta)}{
W^{q}(b-a)}  \right)\frac{W^{(q)}(b-\alpha)}{W^{(q)}(b-a)  }  \right) \right)    \nonumber\\
&=&\frac{W^{(0)}(b-\beta)}{W^{(0)}(b-a)}\int_0^{b-\alpha}W^{0}(x)dx+\frac{W^{(0)}(b-\alpha)}{W^{(0)}(b-a)}\int_{b-\beta}^{b-a}W^{0}(x)dx
-\int_0^{b-\alpha}W^{0}(x)dx  \nonumber \\
&=&K_{a,\alpha,\beta,b}. \label{Def::K}
\end{eqnarray}
Thus, by (\ref{J28}) and Lemma \ref{Prop::Steady} we obtain the desired result.
$\bullet$

\section{An Example}
\label{Sec::Example}

We suppose in this section that
$$Y_t=x +\sigma w_t +N_t,$$
where $N$ is a compound Poisson process independent of $w$ such that the rate of jump arrivals is equal to 1 and the Levy measure $\nu$ of $N$ is
$$\nu(dy)=\theta e^{-\theta y}dy,\quad y\ge 0$$
for some $\theta>0$, and $\nu((-\infty,0])=0$. In this case, we can write the linear operator ${\cal A}$ in the form
\begin{equation}{\cal A}f(x)={\sigma^2\over 2}f''(x) +\int_0^\infty \left[f(x+y)-f(x)\right]\nu(dy).\label{(4)}\end{equation}

We also specify the opportunity cost function as
$$\phi(x)=(x-\rho)^2$$
where $\rho$ is a fixed target value.
Suppose now that $x\in(a,b)$.
 Using (\ref{(4)}), the equation in (i) in Theorem \ref{t2} may be written as
$${1\over 2}\sigma^2f''(x) + (x-\rho)^2-\lambda f(x)+\int_0^\infty \left[f(x+y)-f(x)\right]\theta e^{-\theta y}dy=0.$$
This becomes
$${1\over 2}\sigma^2f''(x) + (x-\rho)^2-(1+\lambda) f(x)+ \theta e^{\theta x}
\int_x^b f(z) e^{-\theta z}dz + $$
\begin{equation}e^{\theta x}\int_b^\infty\left[f(b)+d(z-b)\right]\theta e^{-\theta z}dz=0 \label{(4.1)},
\end{equation}
and also
 \begin{equation}{1\over 2}\sigma^2 e^{-\theta x}   f''(x) + e^{-\theta x}(x-\rho)^2-(1+\lambda) e^{-\theta x}f(x)+ \theta
\int_x^b f(z) e^{-\theta z}dz+\zeta =0,\label{B}\end{equation}
where
$$\zeta = \int_b^\infty\left[f(b)+d(z-b)\right]\theta e^{-\theta z}dz.$$
Let us now introduce $e^{-\theta x}f(x)=g(x)$. We then obtain the following equation from (\ref{B}):
\begin{equation}({1\over 2}\sigma^2\theta^2-\lambda-1)g(x) + \sigma^2\theta g'(x) + {1\over 2}\sigma^2g''(x)+ e^{-\theta x}(x-\rho)^2+\theta\int_x^b g(z)dz + \zeta =0\label{(6)}.\end{equation}
Differentiating the above with respect to $x$ we get the following inhomogeneous linear ordinary differential equation of the third order:
\begin{equation}{1\over 2}\sigma^2g''' + \sigma^2\theta g'' +\left({1\over 2}\sigma^2\theta^2 -\lambda-1\right)g'-\theta g + 2 e^{-\theta x}(x-\rho)-\theta e^{-\theta x}(x-\rho)^2 =0.\label{(7)}\end{equation}
A particular solution for the inhomogeneous equation, denoted by $g_p$, is given by
$$g_p(x)=e^{-\theta x}\left[K_1(x-\rho)^2 + K_2(x-\rho) +K_3\right],$$
where
$$K_1={1\over\lambda}\quad,~ K_2={2\over \theta\lambda^2}\quad,~ K_3={1\over \lambda^3\theta^2}\left[2\lambda+2+\theta^2\lambda\sigma^2\right].$$
The general solution of the homogeneous equation is given by $g_h$, that is,
$$g_h(x)=L_1e^{c_1 x}+ L_2e^{c_2 x}+ L_3e^{c_3 x}$$
where $c_1,c_2,c_3$ are the roots of the equation
$$P(x)={1\over 2}\sigma^2x^3 + \sigma^2\theta x^2+\left({1\over 2}\sigma^2\theta^2 -\lambda-1\right)x-\theta=0$$
and $L_1,L_2,L_3$ are \lq\lq free parameters".
Notice that $P(0)=-\theta<0$, $P(-\theta)=\theta\lambda>0$, and $\lim_{x\to -\infty}P(x)=-\infty$,  $\lim_{x\to \infty}P(x)=\infty$   thus $P(x)$ has three roots, say $c_1<-\theta$, $-\theta<c_2<0$ and $c_3>0$.

We have now arrived at the following family of candidate solutions:
$$g(x;L_1,L_2,L_3,b)=e^{-\theta x}\left[K_1(x-\rho)^2 + K_2(x-\rho) +K_3\right] +L_1e^{c_1 x}+   L_2e^{c_2 x} +
L_3e^{c_3 x}.$$
This gives
\begin{equation}f(x;L_1,L_2,L_3,b)=K_1(x-\rho)^2 + K_2(x-\rho) +K_3 +L_1e^{(\theta+c_1) x}+ L_2e^{(\theta+c_2) x}+
 L_3e^{(\theta+c_3) x}.\label{(7.1)}\end{equation}
For simplicity we shall use the notation $f(x;L_1,L_2,L_3,b)=f(x)$.

We now have 7 unknown parameters $a,\alpha,\beta,b,L_1,L_2, L_3$. From the conditions of Theorem \ref{t2}, we may derive the following 6 equations for these constants:
\begin{equation}f'(a)=-c\label{(8)}\end{equation}
\begin{equation}f'(\alpha)=-c\label{(9)}\end{equation}
\begin{equation}f'(b)=d\label{(10)}\end{equation}
\begin{equation}f'(\beta)=d\label{(11)}\end{equation}
\begin{equation}f(a)=f(\alpha)+C+c(\alpha-a)\label{(12)}\end{equation}
\begin{equation}f(b)=f(\beta)+D+d(b-\beta)\label{(13)}.\end{equation}
In addition, if we trace back our derivation in the above, then we see that we must have (\ref{(4.1)}) hold for at least for one particular $x$ since in going from (\ref{(6)}) to (\ref{(7)}) we took a derivative. Select $x=b$. This then gives us our $7$th equation
\begin{equation}{1\over 2}\sigma^2f''(b) + (b-\rho)^2-(1+\lambda) f(b)+e^{\theta b}\zeta=0.\label{(14)}\end{equation}
We now have the following.

\begin{theorem}\label{t4} Suppose that there exist seven constants $L_1\le 0,L_2\le 0,L_3\le 0, a<\alpha\le\beta<b$ satisfying the seven equations (\ref{(8)})-(\ref{(14)}). We define $h$ by
\begin{equation*}
h(x)=
\begin{cases}
f(a)-c(x-a),&\textrm{if~}x\le a, \\
              f(x),&\textrm{if~} a\le x\le b \\
              f(b)+d(x-b),&\textrm{if~} x\ge b.
\end{cases}
\end{equation*}
Then $h(x)=V(x)$, i.e., $h(x)$ is the value function of the optimization problem.
 Furthermore, the policy $(T^*,\Xi^*)$ described in (\ref{(A2)}) and (\ref{(A3)}) with this choice of $a,\alpha,\beta,b$ is optimal.
\end{theorem}

In order to prove this theorem, we need the following lemma.

\begin{lemma}\label{l5} Assume the conditions of Theorem \ref{t4}. Then there exists a constant $\xi\in(\alpha,\beta)$ such that $h'$ is convex on $[a,\xi]$, concave on $[\xi,b]$. Furthermore $h'(x)\le -c$ if $x\in[a,\alpha]$, $h'(x)\ge d$ if $x\in [\beta,b]$, and $-c\le h'(x)\le d$ if $x\in[\alpha,\beta]$.
\end{lemma}
{\bf Proof:} From the condition that $L_1,L_2,L_3\le 0$ it follows that $h'''(x)$ is decreasing on $(a,b)$. Therefore $h''(x)$ has at most two zero points, which implies that $h'(x)$ has at most two local extreme values in $(a,b)$. By (\ref{(8)})-(\ref{(11)}) the derivative function $h'$ must be    first decreasing then increasing  then again decreasing on $[a,b]$. Since $h''$ is concave, it must be either increasing, or decreasing, or first increasing then decreasing on $[a,b]$. But the first two possibilities are not possible, since $h'$ can not be neither convex nor concave on $[a,b]$. Therefore it must be first convex then concave.

{\bf Proof of Theorem \ref{t4}:} We need to prove that the conditions of Theorem \ref{t2} are satisfied. Condition (i) and the required smoothness of $h$ follows from our construction. Next we prove (ii) and (iv). From conditions (\ref{(9)}), (\ref{(11)}) and Lemma \ref{l5} it follows that
\begin{equation*}
Mh(x)=
\begin{cases}
h(\alpha)+C+c(\alpha-x),& \textrm{if~} a\le x\le\alpha, \\
        h(x)+\hbox{min}\{C,D\},& \textrm{if~} \alpha<x<\beta, \\
         h(\beta)+D+d(x-\beta),& \textrm{if~} \beta\le x\le b.
\end{cases}
\end{equation*}
Condition (ii)  follows from (\ref{(12)}) and (\ref{(13)}).  We show condition (iv) for the case of $x\in[a,\alpha]$, the case of $x\in[\beta,b]$ is similar and the case of $x\in(\alpha,\beta)$ is obvious. We need to show that $0\le h(\alpha)- h(x) +C+c(\alpha-x)$. For $x=a$ we have equality by (\ref{(12)}), and the derivative of the right-hand side of the inequality with respect to $x$ is $-h'(x)-c$, which is non-negative on $[a,\alpha]$ by Lemma
\ref{l5}.
Next we show condition (iii).
First we look at the case of $x>b$. Let
$$K(x)={\cal A}h(x)-\lambda h(x) +(x-\rho)^2;\quad x\in\Re\setminus\{a,b\}.$$
Then
\begin{equation}K(x)={\sigma^2\over 2}h''(x)+\int_0^\infty\left[h(x+y)-h(x)\right]\nu(dy)-\lambda h(x)+(x-\rho)^2=0,\quad x\in(a,b)\label{K1}\end{equation}
and
\begin{equation}K(x)=\int_0^\infty\left[h(x+y)-h(x)\right]\nu(dy)-\lambda h(x)+(x-\rho)^2,\quad x\in(-\infty,a)\cup (b,\infty).\label{K2}\end{equation}
Therefore $0=K(b-)={\sigma^2\over 2}h''(b-)+K(b+)$, and ${\sigma^2\over 2}h''(b-)\le 0$ implies $K(b+)\ge 0$. On the other hand for $x>b$ we have $K'(x)=-\lambda d +2(x-\rho)$ and $K'(b+)=-\lambda d +2(b-\rho)$. We also have by ({\ref{K1}})
\begin{equation}0=K'(x) ={\sigma^2\over 2}h'''(x) +\int_0^\infty \left[h'(x+y)-h'(x)\right]\nu(dy)-\lambda h'(x) +2(x-\rho),\quad x\in(a,b)\label{K13}\end{equation}
hence
$$0=K'(b-)={\sigma^2\over 2}h'''(b-)-\lambda d +2(b-\rho)={\sigma^2\over 2}h'''(b-)+K'(b+).$$
Since $h'''(b-)\le 0$ so we must have $K'(b+)\ge 0$. But $K'$ is increasing on $(b,\infty)$, thus  $K'(x)\ge 0$ whenever $x\in (b,\infty)$. This in turn implies that $K$ is increasing on $(b,\infty)$, thus $K(x)\ge 0$ for $x\in(b,\infty)$.

Next we show that $K(x)\ge 0$ for $x<a$. By (\ref{K1}) and (\ref{K2}) $0=K(a+)=K(a-)+{\sigma^2\over 2}h''(a+)$ which implies that $K(a-)\ge 0$.
Hence all we need to show is that $K'(x)\le 0$ for $x\le a$. Differentiating (\ref{K2}) we get
$$K'(x) =\int_0^\infty \left[h'(x+y)-h'(x)\right]\nu(dy)-\lambda h'(x) +2(x-\rho),\quad x<a,$$
 and with a change of variable in the integral one can see that
\begin{equation}K'(x)=e^{\theta(x-a)}\int_0^\infty[h'(a+z)+c]\nu(dz)+\lambda c+2(x-\rho),\quad x<a\label{K11}\end{equation}
and
\begin{equation}K''(x)=\theta e^{\theta(x-a)}\int_0^\infty[h'(a+z)+c]\nu(dz)+2,\quad x<a.\label{K10}\end{equation}
Thus $K''$ is either increasing or decreasing on $(-\infty,a)$ depending on the sign of $\int_0^\infty[h'(a+z)+c]\nu(dz)$ which makes $K'$ either convex or concave on $(-\infty,a)$. However, the fact that \break $\lim_{x\to -\infty}K'(x)=-\infty$ implies that $K'$ must be concave and $K''$ decreasing on $(-\infty,a)$. In addition, the integral in (\ref{K10}) and in (\ref{K11})  is negative. By (\ref{K13}) and (\ref{K11})
$$0=K'(a+)={\sigma^2\over 2}h'''(a+)+K' (a-)$$
thus $K'(a-)\le 0$ follows from $h'''(a+)\ge 0$. Formula (\ref{K11}) imply that $K'$ has at most one zero-point on $(-\infty,a)$. These facts about $K'$ imply that indeed $K'(x)\le 0$.\hfill$\bullet$

We conclude this section by providing an example showing how the limiting distribution $\pi$ of Proposition \ref{prop:stat} may
be calculated for an arbitrary double bandwidth control policy $(a,\alpha,\beta,b)$. We will assume that $Y_t=\vartheta t + N_t$, where $\vartheta < 0$ and $N_t$ is a compound Poisson process which has jumps at rate one and jump sizes which
are exponentially distributed with rate $\theta$.  Note that by (\ref{density::U}) and Proposition \ref{prop:stat}, it suffices to determine the function  $W^{(0)}= \lim_{q \rightarrow 0}W^{(q)}$. By (\ref{LaplaceW}), the Laplace
transform of $W^{(0)}$ is given by $1/\psi(-s)$ where $\psi(s)$ is the Levy exponent of $Y_t$.   Morever, by (8.1) in \cite{K06} it then follows that
\begin{eqnarray*}
\psi(s) &=& \vartheta s - \int_{0}^{\infty}(1-e^{xs}) \theta e^{-\theta x}dx,
\end{eqnarray*}
for $s < \theta$, which reduces to $\psi(s)=\vartheta s + s(\theta-s)^{-1}$. One may now proceed to verify that
\begin{eqnarray*}
\frac{1}{\psi(-s)} &=& \frac{-\theta}{s(\vartheta s + \vartheta \theta +1)} - \frac{1}{\vartheta s +
\vartheta \theta+1}.
\end{eqnarray*}
In the case in which $\theta \vartheta \neq -1$ inverting each of the terms in the above, one obtains that the function $W^{(0)}$ is given by
\begin{eqnarray*}
W^{(0)}(x) &=& \frac{-\theta}{\theta \vartheta +1} - \left (  \frac{1}{ \vartheta }- \frac{\theta}{\theta \vartheta +1}  \right)
\exp\left(-\frac{\theta \vartheta +1}{\vartheta}x \right).
\end{eqnarray*}
For the case in which $\theta \vartheta =1$, one has that $W^{(0)}(x)=\vartheta^{-1}(\theta x +1)$. Substituting into the formula of Proposition  \ref{prop:stat}, one may now obtain the density of $\pi$.
\section{Acknowledgements}

The authors would like to thank Bert Zwart for his help with Section \ref{Sec::Analysis}.

\appendix

\section{Appenidx}

In the Appendix, we provide a proof of the fact that Ito's rule applies to test functions $f \in{\cal D} $.
For   $f\in{\cal D}$ the second derivative  $f''(x)$ may not exist in points $S=\{x_1,\dots,x_n\}$. We shall call $S$ the set of exceptional points. We extend $f''$ to the entire of $\Re$ by assuming an arbitrary value for $f''(x_i)$. This convention will be used in the rest of this section. The following is then the main result.

\begin{proposition}\label{p6} If $f\in{\cal D}$ and $X$ is a controlled cash on hand process with an arbitrary impulse control $(T,\Xi)=(\tau_1,\tau_2,\dots,\tau_n,\dots,\ \xi_1,\xi_2,\dots,\xi_n,\dots,)$ then Ito's rule holds in its usual form:
$$f(X_t)-f(X_0)=\int_{(0,t]} f'(X_{s-})dX_s + {\sigma^2\over 2}\int_{(0,t]}f''(X_s)ds+$$
\begin{equation}\sum_{0<s\le t}\left\{f(X_s) - f(X_{s-})-f'(X_{s-})\Delta X_s\right\}.\label{(16)}\end{equation}
\end{proposition}

In order to prove Proposition \ref{p6}, we need the following two lemmas.

\begin{lemma}\label{l7} Let $f\in{\cal D}$ with set of  exceptional points $S$. Then there exists a sequence $\left(f_n\right)_{n\ge 1}\subset C^2(\Re)$ such that the following hold;\hfill\break
(i) $f_n(x)\to f(x)$ and $f'_n(x)\to f'(x)$ for every $x\in\Re$ as $n\to\infty$;\hfill\break
(ii) $f''_n(x)\to f''(x)$ for every $x\in\Re\setminus S$ as $n\to\infty$;\hfill\break
(iii) $f_n'$ and $f_n''$ are bounded uniformly in $n$, i.e., $|f_n'(x)|\le C_1$ and $|f_n''(x)|\le C_1$ for some constant $C_1$ and all $n$ and $x\in\Re$.
\end{lemma}

The proof of this lemma can be based on the proof of a similar lemma in {\O}kesendal \cite{OK03}, Appendix D with some obvious modifications.

\begin{lemma}\label{l8}
 If $\sigma\not= 0$ then  for every $x\in\Re$
$$\int_0^\infty 1_{\{x\}}(X_s)ds=0$$
In other words, the Lebesgue measure of the time the controlled cash on hand process spends at level $x$ is zero.
\end{lemma}

{\bf Proof:}
\begin{eqnarray*}
&&E\left[\int_0^\infty 1_{\{x\}}(X_s)ds\right]\\
&=&\int_0^\infty P[X_s=x]ds\\
&\le&\int_0^\infty P[\tau_i=s\ \hbox{for some}\ i]ds + \int_0^\infty P[X_s=x,\ \tau_i\not=s\  \hbox{for all}\ i]ds.
\end{eqnarray*}
We deal with these last two integrals separately.
$$\int_0^\infty P[\tau_i=s\ \hbox{for some}\ i]ds\le \int_0^\infty\sum_{i=1}^\infty P[\tau_i=s]ds=\sum_{i=1}^\infty\int_0^\infty P[\tau_i=s]ds$$
and this last expression is zero because the set $\left\{s\ge 0:\ P[\tau_i=s]>0\right\}$ is either countable or finite.
For the second integral we have
$$\int_0^\infty P[X_s=x,\ \tau_i\not=s\  \hbox{for all}\ i]ds=\sum_{i=0}^\infty\int_0^\infty P\left[X_s=x,\ \tau_i<s<\tau_{i+1}\right]ds.$$
Now it suffices to show that the probability in the right-hand side is zero. Indeed,
\begin{eqnarray*}
&&P\left[X_s=x,\ \tau_i<s<\tau_{i+1}\right]\\
&=&  \int_{[0,s]\times\Re}P\left[X_s=x,\ s<\tau_{i+1}\ |\ \tau_i=u,X_u=y\right]P\left[\tau_i\in du, X_u\in dy\right]\\
&=&\int_{[0,s]\times\Re}P\left[X_s-X_u=x-y,\ s<\tau_{i+1}\ |\ \tau_i=u,X_u=y\right]P\left[\tau_i\in du, X_u\in dy\right]\\
&=&\int_{[0,s]\times\Re}P\left[Y_s-Y_u=x-y,\ s<\tau_{i+1}\ |\ \tau_i=u,X_u=y\right]P\left[\tau_i\in du, X_u\in dy\right]\\
&\le& \int_{[0,s]\times\Re}P\left[Y_s-Y_u=x-y,\ |\ \tau_i=u,X_u=y\right]P\left[\tau_i\in du, X_u\in dy\right]\\
&=&\int_{[0,s]\times\Re}P\left[Y_s-Y_u=x-y\right]P\left[\tau_i\in du, X_u\in dy\right]\\
&=&\int_{[0,s]\times\Re}P\left[Y_s=x\ |\ Y_u=y\right]P\left[\tau_i\in du, X_u\in dy\right]
\end{eqnarray*}
and $P\left[Y_s=x\ |\ Y_u=y\right]=0$ follows from our assumption $\sigma\not= 0$ and Sato \cite{S99}, Theorem 27.4.

We now provide the proof of Proposition \ref{p6}.

{\bf Proof of Proposition \ref{p6}:} Let  $\left(f_n\right)_{n\ge 1}$ be the sequence approximating $f$ in the sense of Lemma \ref{l7}. Ito's rule holds for each $f_n$, i.e.,
\begin{eqnarray}
f_n(X_t)-f_n(X_0)&=&\int_{(0,t]} f_n'(X_{s-})dX_s + {\sigma^2\over 2}\int_{(0,t]}f_n''(X_s)ds \label{(17)}\\
&&+
\sum_{0<s\le t}\left\{f_n(X_s) - f_n(X_{s-})-f_n'(X_{s-})\Delta X_s\right\}. \nonumber
\end{eqnarray}
All we need to show that all three terms in the right-hand side of (\ref{(17)}) converge to the corresponding terms in the right-hand side of (\ref{(16)}) as $n\to\infty$.
We can write $X=X_0+M_1(t)+A_1(t)$ where $M_1$ is a local martingale with bounded jumps (thus also locally square-integrable) and $A_1$ is a finite variation process (Jacod \& Shiryaev \cite{JS03}, Proposition 4.17). We then have
$$\int_{(0,t]} f_n'(X_{s-})dM_1(s) \to \int_{(0,t]} f'(X_{s-})dM_1(s)\ \hbox{in probability as}\ n\to\infty$$
by Theorem 4.40 {iii} in Jacod \& Shiryaev \cite{JS03}. Also
 $$\int_{(0,t]} f_n'(X_{s-})dA_1(s) \to \int_{(0,t]} f'(X_{s-})dA_1(s)\ \hbox{a.s, as}\ n\to\infty$$
by the Dominated Convergence Theorem. Therefore, the first integral in the right-hand side of (\ref{(17)}) indeed converges to the corresponding integral in (\ref{(16)}).
The convergence
$$\sum_{0<s\le t}\left\{f_n(X_s) - f_n(X_{s-})-f_n'(X_{s-})\Delta X_s\right\}\to \sum_{0<s\le t}\left\{f(X_s) - f(X_{s-})-f'(X_{s-})\Delta X_s\right\}$$
follows from the discrete time version of the Dominated Convergence Theorem since $|f_n(X_s) - f_n(X_{s-})-f_n'(X_{s-})\Delta X_s|$ is is bounded by ${C_1\over 2}(\Delta X_s)^2$ and $\sum_{0<s\le t}(\Delta X_s)^2<\infty$.
Finally we need to show that
$${\sigma^2\over 2}\int_{(0,t]}f_n''(X_s)ds \to {\sigma^2\over 2}\int_{(0,t]}f''(X_s)ds$$
as $n\to\infty$. If $\sigma=0$ then there is nothing to prove and if $\sigma\not= 0$ then this follows from Lemma \ref{l8} and the Dominated Convergence Theorem.

\vskip .2in

\bibliographystyle{plain}
\bibliography{references}

\end{document}